\documentclass[11pt]{article}
\usepackage{amsmath,amssymb,amsthm}

\newtheorem{lemma}{Lemma}[section]
\newtheorem{theorem}[lemma]{Theorem}

\theoremstyle{definition}

\newtheorem{remark}[lemma]{Remark}

\begin{document}

\title{Linearization and exponential stability}
\author{Hans Zwart\\Department of Applied Mathematics\\ University of Twente\\P.O.\ Box 217,\\ 7500 AE Enschede, The Netherlands\\
h.j.zwart@utwente.nl}
%\date{April 1, 2014}

\maketitle
\begin{abstract} 
We give sufficient conditions such that the exponential stability of the linearization of a non-linear system implies that the non-linear system is (locally) exponentially stable. One of these conditions is that the non-linear system is Fr\'{e}chet differential at the equilibrium, if it is only Gateaux differentiable, then we show by means of an example that the result does not hold.

\end{abstract}

\section{Introduction}

For finite-dimensional systems it is well-known that if the linearization of a non-linear differential equation around a equilibrium point is exponentially stable, then equilibrium point is locally exponentially stable for the original equation.  We study this question for infinite-dimensional space.  We study the following abstract differential equation
\begin{equation}
  \label{eq:1}
  \dot{x}(t) = A x(t) + f(x(t)), \qquad x(0)=x_0,
\end{equation}
where $A$ is the infinitesimal generator of a $C_0$-semigroup on the
Hilbert space $X$, and $f: X \mapsto X$ is a locally Lipschitz
continuous function with $f(0)=0$.

It is clear that the origin is an equilibrium point of (\ref{eq:1}). In the following  section we study the question whether the exponential stability
of the $C_0$-semigroup generated by $A$ implies the same for the non-linear equation (\ref{eq:1}). In Section \ref{sec:2} we prove a positive result, whereas in Section \ref{sec:3} we show by means of a simple example that the finite dimensional result does not hold in its full generality on infinite-dimension spaces.

\section{Sufficient conditions for exponential stability}
\label{sec:2}

We formulate and prove a positive result. Here we denote the domain of the operator $A$ by $D(A)$, and the class of bounded, linear operators from $X$ to $X$ by ${\mathcal L}(X)$. We say that the semigroup $\left( T(t) \right)_{t \geq 0}$ exponentially stable, when there exists a $M$ and $\omega_0>0$ such that $\|T(t) \| \leq M e^{-\omega_0t}$. 
\begin{theorem}
\label{T:1.1} Let $f$ have zero Fr\'{e}chet derivative at zero. If $A$ generates an exponentially stable semigroup on $X$ and there $\omega \in
{\mathbb R}$ a
bounded, boundedly invertible self-adjoint map $Q \in {\mathcal L}(X)$ and such that for all $x \in D(A)$
\begin{equation}
  \label{eq:2}
  \langle x, QA x \rangle + \langle Ax, Q x \rangle \leq \omega
  \|x\|^2,
\end{equation}
then (\ref{eq:1}) is (locally) exponentially stable around zero.
\end{theorem}
\noindent
{\bf Proof}\/
Without loss of generality we may assume that $\omega > 0$.

Since $A$ generates an exponentially stable semigroup there exists a
self-adjoint, positive $P \in {\mathcal L}(X)$ such that, see \cite[Theorem 5.1.3]{CuZw95},
\begin{equation}
  \label{eq:3}
  \langle x, PA x \rangle + \langle Ax, P x \rangle =- \|x\|^2.
\end{equation}
Now define $P_2= 2 P + \frac{1}{\omega} Q$. This is a bounded and boundedly
invertible self-adjoint linear operator. This implies that there exists $m_1, M_1 >0$ such that for all $x \in X$
\begin{equation}
  \label{eq:6}
  m_1 \|x\|^2 \leq \langle x, P_2 x\rangle \leq M_1 \|x\|^2.
\end{equation}
Combining (\ref{eq:2}) and (\ref{eq:3}) gives that
\begin{equation}
  \label{eq:4}
  \langle x, P_2A x \rangle + \langle Ax, P_2 x \rangle \leq - \|x\|^2.
\end{equation}
Now we choose as Lyapunov function for (\ref{eq:1}) the function
\[
   V(x) = \langle x, P_2 x \rangle.
\]
Using (\ref{eq:4}) we find
\begin{align*}
   \dot{V}(x) =&\ \langle x, P_2(A x+ f(x)) \rangle + \langle Ax+f(x), P_2 x \rangle\\
   \leq &\ - \|x\|^2 + \langle x,P_2 f(x) \rangle +  \langle x,P_2 f(x) \rangle\\
   \leq &\ - \|x\|^2 + 2 \|P_2\| \|x\|\|f(x)\|.
\end{align*}
By (\ref{eq:6}) and since the Fr\'{e}chet derivative of $f$ at zero is zero, we can find a $\delta>0$ such for all $x$ with $V(x) \leq \delta$, there holds that
\begin{equation}
  \label{eq:7}
  \dot{V}(x) \leq \frac{-1}{2M_1} V(x).
\end{equation}
So the abstract differential equation
(\ref{eq:1}) is exponentially stable locally around the origin.\hfill QED
\medskip

\begin{remark}
Condition (\ref{eq:2}) is for many physical examples not a hard
condition. However, since equation (\ref{eq:4}) implies that $A$ is
similar to a dissipative operator, it does not always hold. 

% If I am not mistaken, then mathematicians in numerics has seen the strange behavior that
% the linearization is stable, but the non-linear not. They regard this
% as being far away from normality. The name to mention here is Trefethen.

% The way to build a generator not similar to a contraction semigroup is
% to use a conditional basis.
\end{remark}

In the next section we shall show that if the right hand-side of
(\ref{eq:1}) is just slightly different, then Theorem \ref{T:1.1} does
not hold.

\section{Gateaux linearization exponentially stable, but system not}
\label{sec:3}

In Theorem \ref{T:1.1} we assumed that the Fr\'{e}chet derivative at the origin was zero. By means of an example, we show that this condition cannot be replaced by the condition that the Gateaux derivative at the origin must be zero.

As state space we take $X=\ell^2({\mathbb N})$, and we consider the
differential equation
\begin{equation}
  \label{eq:8}
  \dot{x}(t) = - x(t) + f(x(t)), \qquad x(0)= x_0
\end{equation}
with $f$ given by
\begin{equation}
  \label{eq:9}
  \left(f(x)\right)_n = 3 \sqrt[n]{|x_n|})x_n.
\end{equation}
Hence our system is a diagonal (non-linear) system with on the diagonal
\begin{equation}
  \label{eq:5}
  \dot{x}_n(t) = (-1 + 3 \sqrt[n]{|x_n(t)|})x_n(t). %\qquad n \in
%  {\mathbb N}.
\end{equation}
We summarize results of these scalar differential equations in a lemma. The proof is left to the reader.
\begin{lemma}
\label{L2.1} The differential equation (\ref{eq:5}) has the following properties.
\begin{itemize}
\item The equilibrium's are $\pm 3^{-n}$ and zero.
\item The right hand-side of (\ref{eq:9}) is locally Lipschitz continuous, and for $|x_n|\leq r$ the Lipschitz constant can be majorized by $3(1+\frac{1}{n}) \sqrt[n]{r}$.
\item For $x_n(0) \in (-3^{-n},3^{-n})$ the state converges to zero,
  and for $|x_n(0)| > 3^{-n}$ the state diverges.
\item For $|x_n(0)| > 3^{-n}$ there is a finite escape time.
% \item Note that $\int\left[\frac{1}{-x+3\sqrt[n]{|x|}x}\right] dx = n  \log(|1 - 3 \sqrt[n]{|x|}|) -  \log(|x|)$.
\item The linearization of (\ref{eq:5}) around zero is $\dot{x}_n(t) =
  -x_n(t)$ and thus exponentially stable.
\end{itemize}
\end{lemma}

These result are used to characterize the behavior of the non-linear system (\ref{eq:8}).
\begin{theorem}
\label{T:2.2}
For the non-linear system (\ref{eq:8}) and (\ref{eq:9}) the following holds.
\begin{enumerate}
\item $f$ is (locally) Lipschitz continuous from $X$ to $X$.
\item $f$ is Gateaux differentiable but not Fr\'{e}chet.
\item The origin is an unstable equilibrium point.
\end{enumerate}
\end{theorem}
\noindent{\bf Proof}\/
{\em 1.}\/ Let $x,z$ be two elements of $X$ with norm bounded by $r$. Without loss of generality we may assume that $r>1$. Since the norms are bounded by $r$, the same holds for the absolute value of every element, i.e., $|x_n|, |z_n| \leq r$. Hence we find that
\begin{align*}
  \| f(x)-f(z)\|^2 = &\ \sum_{n=1}^{\infty} | 3 \sqrt[n]{|x_n|})x_n - 3 \sqrt[n]{|z_n|})z_n|^2\\
  \leq&\ \sum_{n=1}^{\infty} |\left( 3 (1+ \frac{1}{n} \sqrt[n]{r}) \right)^2 |x_n - z_n|^2\\
  \leq&\ (6r)^2 \|x-z\|^2,
\end{align*}
where we have used Lemma \ref{L2.1} and the fact that $r>1$. Thus $f$ is Lipschitz continuous, and so is the right hand-side of (\ref{eq:8}).
\smallskip

\noindent
{\em 2..}\/ We show that the Gateaux derivative of $f$ is zero. This implies that the (Gateaux) linearization of (\ref{eq:8}) is $\dot{x}(t) = - x(t)$.

For $x \in X$  and $\varepsilon \in {\mathbb R}\setminus\{0\}$ we have
\begin{eqnarray}
  \label{eq:10}
  \|\frac{f(0+\varepsilon x) - f(0)}{\varepsilon} - 0 \|^2 &=&
  \sum_{n=1}^{\infty} 9 \sqrt[n]{\varepsilon^2 x_n^2} x_n^2\\
\nonumber
  &=& 9 \sum_{n=1}^{\infty} \sqrt[n]{\varepsilon^2} \sqrt[n]{x_n^2} x_n^2
\end{eqnarray}
Next take a $\delta \in (0,1)$ and choose $N$ such that $\sum_{n=N}^{\infty}
x_n(t)^2 \leq \delta$. In particular this implies that
$\sqrt[n]{x_n^2} \leq 1$ for $n\geq N$. Now choose $\varepsilon$ such that
$|\varepsilon|<1$ and $ \sum_{n=1}^{N-1} \sqrt[n]{\varepsilon^2}
\sqrt[n]{x_n^2} x_n^2 \leq \delta$. Combining these two gives that
for this $\varepsilon$ there holds that
\[
   \|\frac{f(0+\varepsilon x) - f(0)}{\varepsilon} - 0 \|^2 \leq 9
   (\delta + \delta).
\]
Since $\delta$ is arbitrarily, this show that
\[
  \lim_{\varepsilon \rightarrow 0} \|\frac{f(0+\varepsilon x) -
    f(0)}{\varepsilon} - 0 \|^2 =0
\]
and so $0$ is the Gateaux derivative of (\ref{eq:5}).

If $f$ would be Fr\'{e}chet differentiable, then its derivative would equal the Gateaux derivative, and thus zero. However, by choosing in equation (\ref{eq:10}) $\varepsilon=1$ and $x= (x_{n})_{n \in {\mathbb N}}$ with $x_{n}=0$ for $n \neq N$ and $x_{N} = 2^{-N}$, we see that $\limsup_{\|x\|\rightarrow 0} \|f(x)\|/\|x\| >0$.
\smallskip

\noindent
{\em 3.}\/
We choose $x(0)= (x_{0n}) _{n \in {\mathbb N}}$ with $x_{0n}=0$ for $n \neq N$ and $x_{0N} = 2^{-N}$. By Lemma \ref{L2.1} we see that the $N$-th equation of (\ref{eq:8}) is unstable, and thus the state $x(t)$ diverge. Since for $N \rightarrow \infty$, there holds $\|x(0)\| \rightarrow 0$,  we see that there exists an initial
state arbitrarily close to zero which is unstable. Thus the non-linear
system is not stable in the origin. \hfill QED
\medskip

The example is this section is not uniformly Lipschitz continuous, and almost every solution of (\ref{eq:8}) will have finite escape time. The following simple adaptation of (\ref{eq:5}) gives a uniformly Lipschitz continuous differential equation on $X$, 
\[
   \dot{x}_n(t) = \frac{ (-1 + 3 \sqrt[n]{|x_n(t)|})x_n(t)}{1 + x_n(t)^2}.
\]

\subsection*{Acknowledgment}
The results in this short note were found after a question by Kirsten Morris. The author would like to thank her and George Weiss for the stimulating discussions.

\end{document}